\documentclass[a4paper, 12pt, leqno]{amsart}
\usepackage{amsmath,amsfonts,amssymb}
\usepackage{mathrsfs}
\usepackage{amscd}
\usepackage{amsthm}
\usepackage{mathrsfs}
\usepackage{hyperref}
\usepackage{epic}
\usepackage{mathtools}

\title{The lowest two-sided cell of weighted Coxeter groups of rank 3}
\author{Jianwei Gao}

\address{Jianwei Gao\\
Beijing International Center for Mathematical Research\\
Peking University\\
No.5 Yiheyuan Road\\
Beijing, 100871\\
People's Republic of China}
\email{gaojianwei@bicmr.pku.edu.cn}

\begin{document}

\maketitle

\begin{abstract}
In this paper, we give precise description for the lowest lowest two-sided cell $c_{0}$ and the left cells in it for a weighted Coxeter group of rank 3. Then we show conjectures $P1-P15$ and $\widetilde{P}$ hold for $c_{0}$ and do some calculation for the based ring of $c_{0}$.
\end{abstract}

\setcounter{section}{-1}
\section{Introduction}

\medskip

Let $(W,S,L)$ be a weighted Coxeter group. In [L3, 13.4], G.Lusztig conjectured that the maximal weight value of the longest elements of the finite parabolic subgroups of $W$ is a bound for $(W,S,L)$. This property is referred as boundness of a weighted Coxeter group([L3, 13.2]). When the rank of $W$ is 1 or 2, this conjecture is clear. In [Gao, 2.1], I proved this conjecture when $rank(W)=3$. As a consequence, $W$ has a lowest two-sided cell $c_{0}$, see [Gao, 7.1]. In this paper, we give precise description for $c_{0}$ and the left cells in it. Then we can show conjectures $P1-P15$ and $\widetilde{P}$ hold for $c_{0}$. At last, we do some calculation for the based ring of $c_{0}$.

\medskip

\section{Preliminaries}

\medskip

\noindent{\bf 1.1.}
In this paper, for any Coxeter group $(W,S)$, we assume the generating set $S$ is finite. We call $|S|$ the rank of $(W,S)$ and denote it by $rank(W)$. We use $l$ for the length function and $\leq$ for the Bruhat order of $W$. The neutral element of $W$ will be denoted by $e$. For $x\in W$, we set $$\mathcal{L}(x)=\{s\in S|sx<x\},\ \mathcal{R}(x)=\{s\in S|xs<x\}.$$
For $s,t \in W$, let $m_{st}\in\mathbb{Z}_{\geq1}\bigcup\{\infty\}$ be the order of $st$ in $W$.

For any $I\subseteq S$, let $W_{I}=\langle I\rangle$. Then $(W_{I},I)$ is also a Coxeter group, called a parabolic subgroup of $(W,S)$. Denote the longest element of $W_{I}$ by $w_{I}$ if $|W_{I}|<\infty$. For $s,t\in S$, $s\neq t$, we use $W_{st}$ instead of $W_{\{s,t\}}$ and $w_{st}$ instead of $w_{\{s,t\}}$.

For $w_{1},w_{2},\cdots,w_{n}\in W$, we often use the notation $w_{1}\cdot w_{2}\cdot\  \cdots \ \cdot w_{n}$ instead of $w_{1}w_{2}\cdots w_{n}$ if
$l(w_{1}w_{2}\cdots w_{n})=l(w_{1})+l(w_{2})+\cdots +l(w_{n})$.

\medskip

\noindent{\bf 1.2.}
Let $(W,S)$ be a Coxeter group. A map $L:W\longrightarrow \mathbb{Z}$ is called a weight function if $L(ww')=L(w)+L(w')$ for any $w,w'\in W$ with $l(ww')=l(w)+l(w')$. Then we call $(W,S,L)$ a weighted Coxeter group. In this paper, the weight function $L$ for any weighted Coxeter group $(W,S,L)$ is assumed to be positive, that is, $L(s)>0$ for any $s\in S$.

\medskip

\noindent{\bf 1.3.}
Let $(W,S,L)$ be a weighted Coxeter group and $\mathbb{Z}[v,v^{-1}]$ be the ring of Laurent polynomials in an indeterminate $v$ with integer coefficients. For $f=\sum \limits_{n\in \mathbb{Z}}a_{n}v^{n}\in \mathbb{Z}[v,v^{-1}]\setminus\{0\}$, we define $deg\ f=\max \limits_{n\in \mathbb{Z}\atop a_{n}\neq 0} n$. Complementally, we define $deg\ 0=-\infty$.

For $w\in W$, set $v_{w}=v^{L(w)}\in \mathbb{Z}[v,v^{-1}]$.
The Hecke algebra $\mathcal{H}$ of $(W,S,L)$ is the unital associative $\mathbb{Z}[v,v^{-1}]-$algebra defined by the generators $T_{s}(s\in S)$ and the relations
$$
(T_{s}-v_{s})(T_{s}+v_{s}^{-1})=0,\ \forall s\in S.
$$
$$
\underbrace{T_{s}T_{t}T_{s}\cdots}_{m_{st}\ factors}\ =\ \underbrace{T_{t}T_{s}T_{t}\cdots}_{m_{st}\ factors},\ \forall s,t\in S,\ m_{st}<\infty.
$$

Obviously, $T_{e}$ is the multiplicative unit of $\mathcal{H}$. For any $w\in W$, we define $T_{w}=T_{s_{1}}T_{s_{2}}\cdots T_{s_{n}}\in \mathcal{H}$, where $w=s_{1}s_{2}\cdots s_{n}$ is a reduced expression of $w$ in $W$. Then $T_{w}$
is independent of the choice of reduced expression and $\{T_{w}|w\in W\}$ is a $\mathbb{Z}[v,v^{-1}]-$basis of $\mathcal{H}$, called the standard basis.
We define $f_{x,y,z}\in \mathbb{Z}[v,v^{-1}]$ for any $x,y,z\in W$ by the identity
$$
T_{x}T_{y}=\sum \limits_{z\in W} f_{x,y,z}T_{z}.
$$
The following involutive automorphism of rings is useful, called the bar involution:
$$\bar{}:\mathcal{H}\longrightarrow \mathcal{H}\ \ \ \ \ $$
$$v^{n}\longrightarrow v^{-n}.$$
$$T_{s}\longrightarrow T_{s}^{-1}.$$
We have $\overline{T_{w}}=T_{w^{-1}}^{-1}$ for any $w\in W$. We set
$$\mathcal{H}_{\leq 0}=\bigoplus \limits_{w\in W}\mathbb{Z}[v^{-1}]T_{w},\ \mathcal{H}_{< 0}=\bigoplus \limits_{w\in W}v^{-1}\mathbb{Z}[v^{-1}]T_{w}.$$
We can get the following facts by easy computation.

\medskip

\noindent{\bf Lemma 1.4.}
(1) For any $x,y,\in W$, we have
$$f_{x,y,e}=\delta_{x,y^{-1}}.$$
(2) For any $x,y,z\in W$, we have
$$deg\ f_{x,y,z}\leq min\{L(x),L(y),L(z)\}.$$
(3) For any $x,y,z\in W$, we have
$$f_{x,y,z^{-1}}=f_{y,z,x^{-1}}=f_{z,x,y^{-1}}.$$
(4) For any finite parabolic subgroup $W_{I}$ of $W$, $x\in W_{I}$, we have
$$deg\ f_{w_{I},w_{I},x}=L(x).$$

\medskip

Define the degree map
$$deg:\ \mathcal{H}\longrightarrow \mathbb{Z}\cup\{-\infty\}\ \ \ \ \ \ \ \ \ $$
$$\ \ \ \ \sum_{w\in W}f_{w}T_{w}\longrightarrow Max\{deg\ f_{w}|w\in W\}.$$
And we set $N=\max \limits_{I\subseteq S\atop |W_{I}|<\infty} L(w_{I})$. G.Lusztig gave the following conjecture in [L3, 13.4].

\medskip

\noindent{\bf Conjecture 1.5.}
Let $(W,S,L)$ be a weighted Coxeter group, $S$ is finite and $L$ is positive, then $N$ is a bound for $(W,S,L)$. Namely, $deg\ (T_{x}T_{y})\leq N$ for all $x,y \in W$.

\medskip

\noindent{\bf Remark 1.6.}
When $W$ is a finite Coxeter group, this conjecture can be proved using Lemma 1.4(2). In [L1, 7.2], G.Lusztig proved this conjecture when $W$ is an affine Weyl group. In [SY, 3.2], J.Shi and G.Yang proved this conjecture when $W$ has complete Coxeter graph. In [Gao, 2.1], I proved this conjecture when the rank of $W$ is 3.

\medskip

\noindent{\bf 1.7.}
For any $w\in W$, there exists a unique element $C_{w}\in \mathcal{H}_{\leq 0}$ such that
$\overline{C_{w}}=C_{w}$ and $C_{w}-T_{w}\in \mathcal{H}_{< 0}$. The elements $\{C_{w}|w\in W\}$ form a $\mathbb{Z}[v,v^{-1}]-$basis of $\mathcal{H}$, called the Kazhdan-Lusztig basis. We define $h_{x,y,z},p_{x,y}\in \mathbb{Z}[v,v^{-1}]$ for any $x,y,z\in W$ such that
$$C_{x}C_{y}=\sum \limits_{z\in W} h_{x,y,z}C_{z}.$$
$$C_{y}=\sum \limits_{x\in W}p_{x,y}T_{x}.$$
These polynomials $p_{x,y}$ are called Kazhdan-Lusztig polynomials. For $w\in W$, $I=\mathcal{R}(w)$, $J=\mathcal{L}(w)$, we have the factorization
$w=x\cdot w_{I}=w_{J}\cdot y$ for some $x,y\in W$. Then we set
$$E_{x}=\sum_{x'\leq x \atop l(x'w_{I})=l(x')+l(w_{I})}p_{x'w_{I},w}T_{x'}.$$
$$F_{y}=\sum_{y'\leq y \atop l(w_{J}y')=l(w_{J})+l(y')}p_{w_{J}y',w}T_{y'}.$$

\medskip

\noindent{\bf 1.8.}
Using Kazhdan-Lusztig basis, we can define the preorders $\underset{L}{\leq}$, $\underset{R}{\leq}$, $\underset{LR}{\leq}$ on $W$. These preorders give rise to equivalence relations $\underset{L}{\sim}$, $\underset{R}{\sim}$, $\underset{LR}{\sim}$ on $W$ respectively. The equicalence classes are called left cells, right cells and two-sided cells of $W$. Then we have partial orders $\underset{L}{\leq}$, $\underset{R}{\leq}$, $\underset{LR}{\leq}$ on the sets of left cells, right cells and two-sided cells of $W$ respectively. For $x,y\in W$, we have $x\underset{L}{\leq}y$ if and only if $x^{-1}\underset{R}{\leq}y^{-1}$.

\medskip

Now we assume Conjecture 1.5 holds and $W$ has a lowest two-sided cell $c_{0}$. It is easy to see that $deg\ h_{x,y,z}\leq N$ for any $x,y,z\in W$, so we can define the a-function
$$a:W\longrightarrow \mathbb{N}\ \ \ \ \ \ \ \ \ \ \ $$
$$\ \ \ \ \ \ \ \ \ \ \ w\longrightarrow \max \limits_{x,y\in W} deg\ h_{x,y,w}.$$

\medskip

For any $x,y,z\in W$, we define $\beta_{x,y,z},\gamma_{x,y,z}\in \mathbb{Z}$ such that
$$f_{x,y,z^{-1}}=\beta_{x,y,z}v^{N}+\mbox{lower degree terms}.$$
$$h_{x,y,z^{-1}}=\gamma_{x,y,z}v^{a(z)}+\mbox{lower degree terms}.$$
\noindent{\bf Lemma 1.9.} Let $x,y,z\in W$.
\\(1) We have $\beta_{x,y,z}=\beta_{y,z,x}=\beta_{z,x,y}$.
\\(2) We have $\gamma_{x,y,z}=\gamma_{y^{-1},x^{-1},z^{-1}}$.
\\(3) If $\beta_{x,y,z}\neq 0$, then $x\underset{L}{\sim}y^{-1}$, $y\underset{L}{\sim}z^{-1}$, $z\underset{L}{\sim}x^{-1}$, $a(x)=a(y)=a(z)=N$, and $\beta_{x,y,z}=\gamma_{x,y,z}=\gamma_{y,z,x}=\gamma_{z,x,y}$.
\\(4) If $\gamma_{x,y,z}\neq 0$ and $a(z)=N$, then $\beta_{x,y,z}=\gamma_{x,y,z}\neq 0$.

\medskip

\noindent{\bf 1.10.}
In [L3, 14.2,14.3], G.Lusztig gave some conjectures about the Hecke algebra, which are called $P1-P15$ and $\widetilde{P}$. When these conjectures hold for $c_{0}$, we can study the based ring $J_{0}$ of $c_{0}$, which is
an associative ring with $\mathbb{Z}-$basis $\{t_{w}|w\in c_{0}\}$ and multiplication
$$t_{x}t_{y}=\sum_{z\in c_{0}}\gamma_{x,y,z^{-1}}t_{z}$$
for any $x,y\in c_{0}$. In section 4 of this paper, we will do some calculation for $J_{0}$ to study its structure when the rank of $W$ is 3.

\medskip

\section{The boundness of $(W,S,L)$}

\medskip

From now on, we assume $(W,S,L)$ is a weighted Coxeter group of rank 3 and $L$ is positive. We set $N=\max \limits_{I\subseteq S\atop |W_{I}|<\infty} L(w_{I})$. In [Gao, 2.1], I proved Conjecture 1.5 in this case.

\medskip

\noindent{\bf Theorem 2.1.}
We have $deg\ (T_{x}T_{y})\leq N$ for all $x,y \in W$.

\medskip

\noindent{\bf Remark 2.2.}
In the proof of [Gao, 2.1], I set $S=\{r,s,t\}$, $m_{rt}=2$, $m_{sr}\geq m_{st}$, $W$ is infinite and $W$ is not an affine Weyl group, otherwise this conjecture has been proved. When $m_{sr}=\infty$ and $m_{st}=2$, or $m_{sr}=m_{st}=\infty$, the conjecture is easy to check, so I only need to consider the following three cases in [Gao].
\\Case 1:  $m_{sr}=\infty>m_{st}\geq 3$.
\\\begin{picture}(60,10)(-45,10)
\put(0,0){\circle{6}}
\put(60,0){\circle{6}}
\put(120,0){\circle{6}}
\put(-4,-15){$r$}
\put(56,-15){$s$}
\put(116,-15){$t$}
\put(3,0){\line(1,0){54}}
\put(63,0){\line(1,0){54}}
\put(25,5){$\infty$}
\put(85,5){$m_{st}$}
\end{picture}
\\\\\\Case 2:  $\infty>m_{sr}\geq m_{st}\geq 4$, $m_{sr}\geq5$.
\\
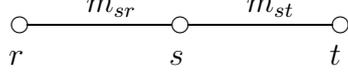
\begin{picture}(60,10)(-45,10)
\put(0,0){\circle{6}}
\put(60,0){\circle{6}}
\put(120,0){\circle{6}}
\put(-4,-15){$r$}
\put(56,-15){$s$}
\put(116,-15){$t$}
\put(3,0){\line(1,0){54}}
\put(63,0){\line(1,0){54}}
\put(25,5){$m_{sr}$}
\put(85,5){$m_{st}$}
\end{picture}
\\\\\\Case 3:  $\infty>m_{sr}\geq7$, $m_{st}=3$.
\\\begin{picture}(60,10)(-45,10)
\put(0,0){\circle{6}}
\put(60,0){\circle{6}}
\put(120,0){\circle{6}}
\put(-4,-15){$r$}
\put(56,-15){$s$}
\put(116,-15){$t$}
\put(3,0){\line(1,0){54}}
\put(63,0){\line(1,0){54}}
\put(25,5){$m_{sr}$}
\put(85,5){$3$}
\end{picture}

\vspace{11 mm}

Now we set
$$M=\{w_{J}|J\subseteq S,\ |W_{J}|<\infty,\ L(w_{J})=N\}.$$
$$\Lambda=\{x\cdot u\cdot y|x,y\in W,\ u\in M\}.$$
For any $w_{J}\in M$, we set
$$B_{J}=\{x\in W|\mathcal{R}(x)\subseteq S\setminus J\}.$$
$$U_{J}=\{y\in W|\mathcal{L}(y)\subseteq S\setminus J\ \mbox{and}\ sw_{J}y\notin \Lambda\ \mbox{for all}\ s\in J\}.$$
Then Theorem 2.1 has the following corollaries.

\medskip

\noindent{\bf Corollary 2.3.}
(1) If $x,y\in W$ satisfy $deg\ (T_{x}T_{y})=N$, then both $x\in \Lambda$ and $y\in \Lambda$.
\\(2) For any $w_{J}\in M$, $q\leq w_{J}$, $x,y\in W$, $\mathcal{R}(x),\mathcal{L}(y)\subseteq S\setminus J$, we have $deg\ (T_{xq}T_{y})\leq N-L(q)$. In particular, $T_{xw_{J}}T_{y}=T_{xw_{J}y}$.
\\(3) For any $w_{J}\in M$, the left cell of $W$ containing $w_{J}$ is $\{x\cdot w_{J}|x\in W\}=\{y\in W|\mathcal{R}(y)=J\}$, the right cell of $W$ containing $w_{J}$ is $\{w_{J}\cdot x|x\in W\}=\{y\in W|\mathcal{L}(y)=J\}$.

\medskip

\noindent{\bf Proof.}
See [Gao, Section 6].\hfill $\square$

\medskip

\noindent{\bf Corollary 2.4.}
Let $w_{J}\in M$, $x\in B_{J}$, $q<w_{J}$.
\\(1) If $y\in B_{J}^{-1}$, $l(q)\geq 2$, then $deg\ (T_{xq}T_{y})<N-L(q)$ for the three cases listed in Remark 2.2.
\\(2) If $y\in U_{J}$, then $deg\ (T_{xq}T_{y})<N-L(q)$ for any weighted Coxeter group of rank 3.

\medskip

\noindent{\bf Proof.}
(1) By the proofs in [Gao, 3.3, 4.7, 5.9], we have the following tables.
\\Case 1
\\\begin{tabular}{|c|c|c|}
\hline
$w_{J}$&$q$&$deg\ (T_{xq}T_{y})$
\\\hline
$w_{st}$&$q<w_{st}$, $l(q)\geq2$&0
\\\hline
\end{tabular}
\\\\Case 2
\\\begin{tabular}{|c|c|c|}
\hline
$w_{J}$&$q$&$deg\ (T_{xq}T_{y})$
\\\hline
$w_{sr}$&$sr$&$\leq Max\{L(st),L(sr)\}$
\\\hline
$w_{sr}$&$rs$&$\leq Max\{L(st),L(sr)\}$
\\\hline
$w_{sr}$&$rsr$&$\leq L(s)$
\\\hline
$w_{sr}$&$srs$&$\leq Max\{L(t),L(r)\}$
\\\hline
$w_{sr}$&$q<w_{sr}$, $l(q)\geq4$&0
\\\hline
$w_{st}$&$st$&$\leq L(sr)$
\\\hline
$w_{st}$&$ts$&$\leq L(sr)$
\\\hline
$w_{st}$&$tst$&0
\\\hline
$w_{st}$&$sts$&$\leq L(r)$
\\\hline
$w_{st}$&$q<w_{st}$, $l(q)\geq4$&0
\\\hline
\end{tabular}
\\\\Case 3
\\\begin{tabular}{|c|c|c|}
\hline
$w_{J}$&$q$&$deg\ (T_{xq}T_{y})$
\\\hline
$w_{sr}$&$sr$&$\leq L(srsr)$
\\\hline
$w_{sr}$&$rs$&$\leq L(srsr)$
\\\hline
$w_{sr}$&$srs$&$\leq L(srs)$
\\\hline
$w_{sr}$&$rsr$&$\leq L(srs)$
\\\hline
$w_{sr}$&$srsr$&$\leq L(sr)$
\\\hline
$w_{sr}$&$rsrs$&$\leq L(sr)$
\\\hline
$w_{sr}$&$rsrsr$&$\leq L(s)$
\\\hline
$w_{sr}$&$srsrs$&$0$
\\\hline
$w_{sr}$&$q<w_{sr}$, $l(q)\geq 6$&0
\\\hline
\end{tabular}
\\\\We have $deg\ (T_{xq}T_{y})<N-L(q)$ in all these three cases.
\\(2) It is obvious when $W$ is a finite Coxeter group. We can get this inequality by [Xie1, 3.1] when $W$ is an affine Weyl group and by [Xie2, 3.5] when $W$ has complete Coxeter graph. Now we only consider other case.

If $l(q)=0$, by Corollary 2.3(1), we get $deg\ (T_{x}T_{y})<N$ since $y\notin\Lambda$. If $l(q)=1$, by Corollary 2.3(1), we get $deg\ (T_{xq}T_{qy})<N$ since $qy\notin\Lambda$, so $deg\ (T_{xq}T_{y})<N-L(q)$.
If $l(q)\geq2$, then we must in the three cases listed in Remark 2.2, we get $deg\ (T_{xq}T_{y})<N-L(q)$ by (1).\hfill $\square$

\medskip

\section{The lowest two-side cell $c_{0}$}

\medskip

Now we fix an element $w_{J}\in M$ and let $c_{0}$ be the two-sided cell of $W$ containing $w_{J}$. In [Gao, 7.1], I proved the following consequence.

\medskip

\noindent{\bf Proposition 3.1.}
(1) The two-sided cell $c_{0}$ is the lowest two-sided cell of $W$.
\\(2) We have $\Lambda=\{w\in W|a(w)=N\}\subseteq c_{0}$.

\medskip

By Corollary 2.3 and Corollary 2.4, we can prove the following proposition using the same method used in [Xie2, 3.6,3.7].

\medskip

\noindent{\bf Proposition 3.2.}
(1) For $w_{J}\in M$, $x\in B_{J}$, $y\in U_{J}$, we have $C_{xw_{J}y}=E_{x}C_{w_{J}}F_{y}$.
\\(2) For $w_{J}\in M$, $y\in U_{J}$, $x\in U_{J}^{-1}$, the set $\Gamma_{J,y}=B_{J}w_{J}y$ is a left cell of $W$ and $\Phi_{x,J}=xw_{J}B_{J}^{-1}$ is a right cell of $W$.

\medskip

Now we can give precise description for $c_{0}$ and the left cells in it.

\medskip

\noindent{\bf Theorem 3.3.}
(1) We have $c_{0}=\Lambda=\{w\in W|a(w)=N\}$.
\\(2) The lowest two-sided cell $c_{0}$ can be decomposed into left cells as
$$c_{0}=\coprod_{w_{J}\in M,y\in U_{J}}\ \Gamma_{J,y}\ .$$

\medskip

\noindent{\bf Proof.}
(1) It is clear that $\Lambda=\bigcup\limits_{w_{J}\in M}B_{J}w_{J}U_{J}=\bigcup\limits_{w_{J}\in M}U_{J}^{-1}w_{J}B_{J}^{-1}\subseteq c_{0}$.

We claim that for any $x\in W$, $y\in\Lambda$, if $x\underset{L}{\leftarrow}y$ or $x\underset{R}{\leftarrow}y$, then $x\in\Lambda$. We only prove the $x\underset{L}{\leftarrow}y$ case.
First, we assume $y=z_{1}w_{J}z_{2}$ for some $w_{J}\in M$, $z_{1}\in B_{J}$, $z_{2}\in U_{J}$, $h_{s,y,x}\neq 0$ for some $s\in S$. By Proposition 3.2(1), we have
$$
\begin{aligned}
&\ C_{s}C_{z_{1}w_{J}z_{2}}\\
=&\ C_{s}C_{z_{2}w_{J}}F_{z_{2}}\\
=&\ \sum_{z_{3}\in B_{J}}h_{s,z_{1}w_{J},z_{3}w_{J}}C_{z_{3}w_{J}}F_{z_{2}}\\
=&\ \sum_{z_{3}\in B_{J}}h_{s,z_{1}w_{J},z_{3}w_{J}}C_{z_{3}w_{J}z_{2}}.
\end{aligned}
$$
We get $x=z_{3}w_{J}z_{2}\subseteq \Lambda$ for some $z_{3}\in B_{J}$ and $z_{2}\in U_{J}$, the claim is proved.

For $w\in c_{0}$, we have $w\underset{LR}{\leq}w_{J}$ and $w_{J}\in\Lambda$, so $w\in\Lambda$ by the claim. On the other hand, we have $\Lambda=\{w\in W|a(w)=N\}\subseteq c_{0}$ by Proposition 3.1(2). This conclusion follows.
\\(2) It is clear that $c_{0}=\bigcup\limits_{w_{J}\in M,y\in U_{J}}\ \Gamma_{J,y}\ .$ If $w_{J},w_{J'}\in M$, $y\in U_{J}$, $y'\in U_{J'}$ such that $\Gamma_{J,y}=\Gamma_{J',y'}$,
then there exists $x\in B_{J}$ such that $xw_{J}y=w_{J'}y'$. If $x\neq e$, we take $s\in \mathcal{L}(x)$, then
$s\in \mathcal{L}(xw_{J}y)=\mathcal{L}(w_{J'}y')=J'$ and $sw_{J'}y'=sxw_{J}y\in \Lambda$, contradict to $y'\in U_{J'}$. So we have $x=e$.
Thus, $J=\mathcal{L}(w_{J}y)=\mathcal{L}(w_{J'}y')=J'$, $y=y'$. Therefore, it is a disjoint union.\hfill $\square$

\medskip

\noindent{\bf Proposition 3.4.}
We have the following table for weighted Coxeter groups of rank 3.
\\\\\begin{tabular}{|c|c|}
\hline
$W$&the number of left cells in $c_{0}$
\\\hline
finite Coxeter group&1
\\\hline
affine Weyl group&$|W_{0}|$
\\\hline
$m_{sr}=\infty$, $m_{st}=m_{rt}=2$&2
\\\hline
$m_{sr}=m_{st}=\infty$, $\infty>m_{rt}\geq 2$&$\infty$ or $2m_{rt}$
\\\hline
$m_{sr}=m_{st}=m_{rt}=\infty$&$\infty$ or 3 or 4
\\\hline
other cases&$\infty$
\\\hline
\end{tabular}
\\\\When $W$ is an affine Weyl group, $W_{0}$ denotes the finite Weyl group corresponding to $W$. When the number of left cells in $c_{0}$ has various possibilities, it depends on the weight function $L$.

\medskip

\noindent{\bf Proof.}
It is clear when $W$ is a finite Coxeter group. See [Xie1, 3.7] when $W$ is an affine Weyl group. Then we assume $\infty\geq m_{sr}\geq m_{st}\geq m_{rt}\geq 2$ and consider the following 11 cases.
\\(1) $m_{sr}=m_{st}=m_{rt}=\infty$.
We have $N=Max\{L(r),L(s),L(t)\}$. We may assume $L(r)\geq L(s)\geq L(t)$. If $L(r)=L(s)=L(t)$, then $\Gamma_{\{r\},e}$,$\Gamma_{\{s\},e}$ and $\Gamma_{\{t\},e}$ are all the 3 left cells in $c_{0}$. If $L(r)=L(s)>L(t)$, then $\Gamma_{\{r\},e}$,$\Gamma_{\{s\},e}$,$\Gamma_{\{r\},t}$ and $\Gamma_{\{s\},t}$ are all the 4 left cells in $c_{0}$. If $L(r)>L(s)\geq L(t)$, then for different $k\in \mathbb{N}$, $\Gamma_{\{r\},(st)^{k}}$ are different left cells in $c_{0}$.
\\(2) $m_{sr}=m_{st}=\infty$, $\infty>m_{rt}\geq 2$.
We have $N=Max\{L(w_{rt}),L(s)\}$. If $L(w_{rt})>L(s)$, then for different $k\in \mathbb{N}$, $\Gamma_{\{r,t\},(sr)^{k}}$ are different left cells in $c_{0}$. If $L(w_{rt})=L(s)$, then $\Gamma_{\{r,t\},e}$ and $\Gamma_{\{s\},w},w\in W_{rt}\setminus \{w_{rt}\}$ are all the $2m_{rt}$ left cells in $c_{0}$. If $L(s)>L(w_{rt})$, then $\Gamma_{\{s\},w},w\in W_{rt}$ are all the $2m_{rt}$ left cells in $c_{0}$.
\\(3) $m_{sr}=\infty$, $m_{st}=m_{rt}=2$.
Now we have $N=Max\{L(st),L(rt)\}$. If $L(s)>L(r)$, then $\Gamma_{\{s,t\},e}$ and $\Gamma_{\{s,t\},r}$ are all the 2 left cells in $c_{0}$. If $L(s)=L(r)$, then $\Gamma_{\{s,t\},e}$ and $\Gamma_{\{r,t\},e}$ are all the 2 left cells in $c_{0}$. If $L(r)>L(s)$, then $\Gamma_{\{r,t\},e}$ and $\Gamma_{\{r,t\},s}$ are all the 2 left cells in $c_{0}$.
\\(4) $m_{sr}=\infty$, $\infty>m_{st}\geq3$, $m_{rt}=2$.
We have $N=Max\{L(w_{st}),L(rt)\}$. If $L(w_{st})>L(rt)$, then for different $k\in \mathbb{N}$, $\Gamma_{\{s,t\},(rs)^{k}}$ are different left cells in $c_{0}$. If $L(rt)\geq L(w_{st})$, then for different $k\in \mathbb{N}$, $\Gamma_{\{r,t\},(sr)^{k}}$ are different left cells in $c_{0}$.
\\(5) $m_{sr}=\infty$, $\infty>m_{st},m_{rt}\geq3$.
We have $N=Max\{L(w_{st}),L(w_{rt})\}$. We may assume $N=L(w_{st})$, then for different $k\in \mathbb{N}$, $\Gamma_{\{s,t\},(rst)^{k}}$ are different left cells in $c_{0}$.
\\(6) $\infty>m_{sr}\geq m_{st}\geq4$, $m_{sr}\geq5$, $m_{rt}=2$.
We have $N=Max\{L(w_{sr}),L(w_{st})\}>L(rt)$. If $L(w_{sr})\geq L(w_{st})$, then for different $k\in \mathbb{N}$, $\Gamma_{\{s,r\},(tsrs)^{k}}$ are different left cells in $c_{0}$. If $L(w_{st})>L(w_{sr})$, then for different $k\in \mathbb{N}$, $\Gamma_{\{s,t\},(rst)^{k}}$ are different left cells in $c_{0}$.
\\(7) $\infty>m_{sr}\geq7$, $m_{st}=3$, $m_{rt}=2$.
We have $N=L(w_{sr})>Max\{L(w_{st}),L(rt)\}$. For different $k\in \mathbb{N}$, $\Gamma_{\{s,r\},(tsrsrs)^{k}}$ are different left cells in $c_{0}$.
\\(8) $\infty>m_{sr}\geq m_{st}\geq m_{rt}\geq 4$. If $L(w_{sr})=N$, then for different $k\in \mathbb{N}$, $\Gamma_{\{s,r\},(tsr)^{k}}$ are different left cells in $c_{0}$. $L(w_{st})=N$ and $L(w_{rt})=N$ are similar.
\\(9) $\infty>m_{sr}\geq m_{st}\geq 4$, $m_{rt}=3$, $L(w_{rt})=N$. For different $k\in \mathbb{N}$, $\Gamma_{\{r,t\},(srt)^{k}}$ are different left cells in $c_{0}$.
\\(10) $\infty>m_{sr}\geq m_{st}\geq 4$, $m_{rt}=3$, $L(w_{sr})=N$. For different $k\in \mathbb{N}$, $\Gamma_{\{s,r\},(tsr)^{k}}$ are different left cells in $c_{0}$.
\\(11) $\infty>m_{sr}\geq 4$, $m_{st}=m_{rt}=3$. Now we have $N=L(w_{sr})>L(w_{st})=L(w_{rt})$. For different $k\in \mathbb{N}$, $\Gamma_{\{s,r\},(tsr)^{k}}$ are different left cells in $c_{0}$.
\hfill $\square$

\medskip

For $n\in \mathbb{Z}$, define
$$\pi_{n}:\mathbb{Z}[v,v^{-1}]\longrightarrow \mathbb{Z}\ \ \ \ \ \ \ \ $$
$$\sum_{k\in \mathbb{Z}}a_{k}v^{k}\longrightarrow a_{n}.$$
For $z\in W$, define $\Delta(z)\in \mathbb{N}$ and $n_{z}\in \mathbb{Z}\setminus\{0\}$, such that
$$p_{e,z}=n_{z}v^{-\Delta(z)}+\mbox{lower degree terms}.$$
Let
$$\mathcal{D}=\{z\in W|a(z)=\Delta(z)\}.$$
\\\noindent{\bf Proposition 3.5.}
We have the following propositions $P1'-P15'$ and $\widetilde{P}'$. In particular, Conjectures $P1-P15$ and $\widetilde{P}$ hold for $c_{0}$.
\\$P1'$  For $z\in c_{0}$, we have $a(z)\leq\Delta(z)$.
\\$P2'$  For $d\in \mathcal{D}\cap c_{0}$, $x,y\in c_{0}$, if $\gamma_{x,y,d}\neq0$, then $x=y^{-1}$.
\\$P3'$  For $y\in c_{0}$, there exists unique $d\in \mathcal{D}\cap c_{0}$ such that $\gamma_{y^{-1},y,d}\neq0$.
\\$P4'$  For $z,z'\in W$, if $z\in c_{0}$ or $z'\in c_{0}$, and $z'\underset{LR}{\leq}z$, then $a(z')\geq a(z)$.
\\$P5'$  For $d\in \mathcal{D}\cap c_{0}$, $y\in c_{0}$, if $\gamma_{y^{-1},y,d}\neq0$, then $\gamma_{y^{-1},y,d}=n_{d}=1$.
\\$P6'$  For $d\in \mathcal{D}\cap c_{0}$, we have $d^{2}=e$.
\\$P7'$  For $x,y,z\in W$, if $x\in c_{0}$ or $z\in c_{0}$, then $\gamma_{x,y,z}=\gamma_{y,z,x}$.
\\$P8'$  For $x,y\in W$, $z\in c_{0}$, if $\gamma_{x,y,z}\neq 0$, then $x\underset{L}{\sim}y^{-1}$, $y\underset{L}{\sim}z^{-1}$, $z\underset{L}{\sim}x^{-1}$.
\\$P9'$  For $z,z'\in W$, if $z\in c_{0}$ or $z'\in c_{0}$, $z'\underset{L}{\leq}z$, and $a(z')=a(z)$, then $z'\underset{L}{\sim}z$.
\\$P10'$  For $z,z'\in W$, if $z\in c_{0}$ or $z'\in c_{0}$, $z'\underset{R}{\leq}z$, and $a(z')=a(z)$, then $z'\underset{R}{\sim}z$.
\\$P11'$  For $z,z'\in W$, if $z\in c_{0}$ or $z'\in c_{0}$, $z'\underset{LR}{\leq}z$, and $a(z')=a(z)$, then $z'\underset{LR}{\sim}z$.
\\$P12'$  If $y\in c_{0}\cap W_{I}$ for some $I\subseteq S$, then $a(y)$ computed in terms of $W_{I}$ is equal to $a(y)$ computed in terms of $W$.
\\$P13'$  For any left cell $\Gamma$ in $c_{0}$, we have $|\Gamma\cap\mathcal{D}|=1$. Assume $\Gamma\cap\mathcal{D}=\{d\}$, then $\gamma_{x^{-1},x,d}\neq0$ for all $x\in \Gamma$.
\\$P14'$  For any $z\in c_{0}$, we have $z\underset{LR}{\sim}z^{-1}$.
\\$P15'$  For $x,x'\in W$, $y,w\in c_{0}$, we have the equality $\mathop{\sum}\limits_{y'\in c_{0}}h_{x,y',y}\otimes h_{w,x',y'}=\mathop{\sum}\limits_{y'\in c_{0}}h_{x,w,y'}\otimes h_{y',x',y}$ in $\mathbb{Z}[v,v^{-1}]\mathop{\otimes}\limits_{\mathbb{Z}}\mathbb{Z}[v,v^{-1}]$.
\\$\widetilde{P}'$  For $x,y,z'\in W$, $z\in c_{0}$, if $\gamma_{x,y,z^{-1}}\neq0$, $z'\underset{L}{\leftarrow}z$, then there exists $x'\in c_{0}$, such that $\pi_{a(z)}(h_{x',y,z'})\neq0$.

\medskip

\noindent{\bf Proof.}
By Proposition 3.2 and Theorem 3.3, we can prove $P1'-P15'$ using the same method used in [Xie1, 4.1]. Here we only give the proof of $\widetilde{P}'$.

Since $\gamma_{x,y,z^{-1}}\neq0$, $z\in c_{0}$, by $P8'$, we have $y\underset{L}{\sim}z$. Since $z'\underset{L}{\leftarrow}z$, $z\in c_{0}$, we get $z'\underset{L}{\sim}z$ by $P9'$. So $y,z,z'$ are in the same left cell in $c_{0}$. By Theorem 3.3(2), we assume $y=aw_{J}b$, $z=a'w_{J}b$, $z'=a''w_{J}b$ for some $w_{J}\in M$, $b\in U_{J}$, $a,a',a''\in B_{J}$.
Then we have $deg\ (T_{y}T_{z'^{-1}})=deg\ (T_{aw_{J}b}T_{b^{-1}w_{J}a''^{-1}})=N$. Thus there exists $x'\in W$, such that $deg\ f_{y,z'^{-1},x'^{-1}}=N$, By Lemma 1.9(3) and Theorem 3.3(1), we get $x'\in c_{0}$.
We have $deg\ f_{x',y,z'}=deg\ f_{y,z'^{-1},x'^{-1}}=N$, so $deg\ h_{x',y,z'}=N$, that is, $\pi_{a(z)}(h_{x',y,z'})\neq0$.\hfill $\square$

\medskip

\section{The based ring of $c_{0}$}

\medskip

In this section, for a weighted Coxeter group of rank 3, we study the based ring $J_{0}$ of its lowest two-sided cell $c_{0}$. We want to calculate $t_{x}t_{y}$ for any $x,y\in c_{0}$, but it is quite a difficult problem. For $w_{J},w_{J'}\in M$, we set
$$P_{J,J'}=\{x\in W|\mathcal{L}(x)=J,\ \mathcal{R}(x)=J'\}.$$
Let
$$P=\bigcup_{w_{J},w_{J'}\in M}P_{J,J'}.$$
By [Xie2, 3.10,3.12], we have the following lemma.

\medskip

\noindent{\bf Lemma 4.1.}
(1) For any $w\in c_{0}$, there exists unique $w_{J},w_{J'}\in M$, $p_{w}\in P_{J,J'}$, $x\in U_{J}^{-1}$, $y\in U_{J'}$, such that $w=x\cdot p_{w}\cdot y$. Moreover, we have $C_{w}=E_{x}C_{p_{w}}F_{y}$.
\\(2) Let $w_{1},w_{2},w_{3}\in c_{0}$. We have the following decomposition as in (1).
$$w_{1}=x_{1}\cdot p_{1}\cdot y_{1},\ w_{2}=x_{2}\cdot p_{2}\cdot y_{2},\ w_{3}=x_{3}\cdot p_{3}\cdot y_{3}.$$
If $y_{1}=x_{2}^{-1}$, $y_{2}=x_{3}^{-1}$, $y_{3}=x_{1}^{-1}$, $\mathcal{R}(p_{1})=\mathcal{L}(p_{2})$, $\mathcal{R}(p_{2})=\mathcal{L}(p_{3})$, $\mathcal{R}(p_{3})=\mathcal{L}(p_{1})$, then $\gamma_{w_{1},w_{2},w_{3}}=\gamma_{p_{1},p_{2},p_{3}}$. Otherwise, $\gamma_{w_{1},w_{2},w_{3}}=0$.

\medskip

The lemma above simplify our work to calculate $t_{x}t_{y}$ for $x\in P_{J,J'}$, $y\in P_{J',J''}$, but it is still not easy. So we assume $x$ indecomposable at first.
For $x\in P_{J,J'}\subseteq P$, we call $x$ indecomposable, if $x\notin M$, and there doesn't exist $w_{J''}\in M$, $x_{1}\in P_{J,J''}\setminus M$, $x_{2}\in P_{J'',J'}\setminus M$ such that $x=x_{1}w_{J''}x_{2}$.
Then we have the following proposition.

\medskip

\noindent{\bf Proposition 4.2.}
Let $(W,S,L)$ be a weight Coxeter group of rank 3 but not an affine Weyl group. Assume $w_{J},w_{J'},w_{J''}\in M$, $x\in P_{J,J'}$ and indecomposable, $y\in P_{J',J''}$, then $t_{x}t_{y}=t_{xw_{J'}y}+\delta t_{w_{J}xy}$. If there exists $y_{1}\in W$ such that $y=x^{-1}\cdot y_{1}$, then $\delta=1$, otherwise $\delta=0$.

\medskip

To prove this proposition, we need two lemmas.

\medskip

\noindent{\bf Lemma 4.3.}
Let $x=x_{1}\cdot w_{J'}\in P_{J,J'}$ and indecomposable, $r\in J'$. Then we have $x_{1}\cdot r\notin c_{0}$ or $x_{1}\cdot r\in w_{J}U_{J}$.

\medskip

\noindent{\bf Proof.}
We assume $x_{1}\cdot r\in c_{0}$. Since $x_{1}\cdot w_{J'}\underset{R}{\leq} x_{1}\cdot r$, we get $x_{1}\cdot w_{J'}\underset{R}{\sim} x_{1}\cdot r$ by $P10'$,
so $\mathcal{L}(x_{1}\cdot r)=\mathcal{L}(x_{1}\cdot w_{J'})=J$. Assume $x_{1}\cdot r=w_{J}\cdot x_{2}$ for some $x_{2}\in B_{J}^{-1}$. We claim that $x_{2}\in U_{J}$.
Otherwise, there exists $r'\in J$ such that $r'w_{J}\cdot x_{2}\in c_{0}$.
By Theorem 3.3, we may assume $r'w_{J}\cdot x_{2}=x_{3}\cdot w_{J''}\cdot x_{4}$ for some $w_{J''}\in M$, $x_{3}\in B_{J''}$, $x_{4}\in B_{J''}^{-1}$. Then we have
$$
\begin{aligned}
x=x_{1}\cdot w_{J'}&=x_{1}\cdot r\cdot rw_{J'}\\
&=r'\cdot r'w_{J}\cdot x_{2}\cdot rw_{J'}\\
&=r'\cdot x_{3}\cdot w_{J''}\cdot x_{4}\cdot rw_{J'}\\
&=(r'\cdot w_{3}\cdot w_{J''})w_{J''}(w_{J''}\cdot x_{4}\cdot rw_{J'}).
\end{aligned}
$$
Since $x_{1}\cdot r \underset{R}{\leq} r'\cdot x_{3}\cdot w_{J''}$, we get $x_{1}\cdot r \underset{R}{\sim} r'\cdot x_{3}\cdot w_{J''}$ by $P10'$, so $\mathcal{L}(r'\cdot x_{3}\cdot w_{J''})=\mathcal{L}(x_{1}\cdot r)=J$.
Thus, we get $r'\cdot x_{3}\cdot w_{J''}\in P_{J,J''}$.

\medskip

On the other hand, since $x_{1}\cdot r \underset{L}{\leq} w_{J''}\cdot x_{4}$, by $P9'$, we get $r\in \mathcal{R}(w_{J''}\cdot x_{4})=\mathcal{R}(x_{1}\cdot r)$,
so $\mathcal{R}(w_{J''}\cdot x_{4}\cdot rw_{J'})=J'$. Thus, $w_{J''}\cdot x_{4}\cdot rw_{J'}\in P_{J'',J'}$, contradict to the indecomposability of $x$.\hfill $\square$

\medskip

\noindent{\bf Lemma 4.4.}
Let $x\in P_{J,J'}$ with $J'=\{r_{1},r_{2}\}$ and $r_{1}\neq r_{2}$. Assume $x=x_{1}\cdot w_{J'}$ for some $x_{1}\in B_{J'}$. If $\mathcal{L}(x_{1})\neq J$, then we have $\mathcal{L}(x_{1}r_{1})\neq J$ or $\mathcal{L}(x_{1}r_{2})\neq J$.

\medskip

\noindent{\bf Proof.}
If $\mathcal{L}(x_{1}r_{1})=\mathcal{L}(x_{1}r_{2})=J$, we assume
$$x_{1}\cdot r_{1}=w_{J}\cdot x_{2} \mbox{ for some }x_{2}\in B_{J}^{-1}.$$
$$x_{1}\cdot r_{2}=w_{J}\cdot x_{3} \mbox{ for some }x_{3}\in B_{J}^{-1}.$$
Then we have
$$x=w_{J}\cdot x_{2}\cdot r_{1}w_{J'}=w_{J}\cdot x_{3}\cdot r_{2}w_{J'}.$$
So
$$x_{2}\cdot r_{1}w_{J'}=x_{3}\cdot r_{2}w_{J'}.$$
Since $\mathcal{L}(x_{1})\neq J$, we get $x_{2},x_{3}\in B_{J'}$. Thus, $r_{1}w_{J'}=r_{2}w_{J'}$ and then $r_{1}=r_{2}$, a contradiction.\hfill $\square$

\medskip

\noindent{\bf 4.5.}
Now we prove Proposition 4.2. Finite Coxeter groups don't have any indecomposable elements. If $W$ has complete Coxeter graph, see [Xie2, 4.2]. For $x,y,z\in c_{0}$, we have $\gamma_{x,y,z}=\beta_{x,y,z}$, so we may compute $T_{x}T_{y}$. It is much easier than computing $C_{x}C_{y}$.
\\If $m_{sr}=\infty$ and $m_{st}=m_{rt}=2$, we can list all the indecomposable elements and check this proposition directly.
\\\\\begin{tabular}{|c|c|}
\hline
weight function&indecomposable elements
\\\hline
$L(s)>L(r)$&$srst$
\\\hline
$L(s)=L(r)$&$srt, rst$
\\\hline
$L(r)>L(s)$&$rsrt$
\\\hline
\end{tabular}
\\\\The case of $m_{sr}=m_{st}=\infty$, $m_{rt}=2$ is also easy.
\\\\\begin{tabular}{|c|c|}
\hline
weight function&indecomposable elements
\\\hline
$L(rt)>L(s)$&$rtsrt$
\\\hline
$L(rt)=L(s)$&$rts, srt, srs, sts$
\\\hline
$L(s)>L(rt)$&$srs, sts, srts$
\\\hline
\end{tabular}
\\\\From now on, we only consider the three cases listed in Remark 2.2. We assume $x=x'\cdot w_{J'}$ for some $x'\in B_{J'}$, $y=w_{J'}\cdot y'$ for some $y'\in B_{J'}^{-1}$. Then we have
$$T_{x}T_{y}=\sum_{q\leq w_{J'}}f_{w_{J'},w_{J'},q}T_{x'q}T_{y'}.$$
When $q=w_{J'}$, by Lemma 1.4(4) and Corollary 2.3(2), we have
$$deg\ f_{w_{J'},w_{J'},w_{J'}}=L(w_{J'})=N \mbox{ and } T_{x'q}T_{y'}=T_{xw_{J'}y}.$$
Since $f_{w_{J'},w_{J'},w_{J'}}$ is a monic polynomial, we know
$$f_{w_{J'},w_{J'},w_{J'}}T_{x'w_{J'}}T_{y'}$$
contributes $t_{xw_{J'}y}$ for $t_{x}t_{y}$.
\\When $q<w_{J'}$ and $l(q)\geq 2$, by Corollary 2.4(1), we have $$deg\ (f_{w_{J'},w_{J'},q}T_{x'q}T_{y'})< N.$$
So
$$\sum_{q<w_{J'},\ l(q)\geq 2}f_{w_{J'},w_{J'},q}T_{x'q}T_{y'}$$
has no contribution for $t_{x}t_{y}$.
\\Now we assume $J'=\{r_{1},r_{2}\}\subseteq S$ with $r_{1}\neq r_{2}$. Then we have
$$
\begin{aligned}
&\ \ \ \ \sum_{l(q)\leq 1}f_{w_{J'},w_{J'},q}T_{x'q}T_{y'}\\
&=(v_{r_{1}}-v_{r_{1}}^{-1})T_{x'r_{1}}T_{y'}+(v_{r_{2}}-v_{r_{2}}^{-1})T_{x'r_{2}}T_{y'}+T_{x'}T_{y'}.
\end{aligned}
$$

If $deg\ ((v_{r_{1}}-v_{r_{1}}^{-1})T_{x'r_{1}}T_{y'})=N$ or $deg\ (T_{x'}T_{y'})=N$, then we have $deg\ (T_{x'r_{1}}T_{r_{1}y'})=N$, so $x'r_{1}\in c_{0}$ by Corollary 2.3(1).
Moreover, we have $x'r_{1}\in w_{J}U_{J}$ by Lemma 4.3. Assume $f_{x'r_{1},r_{1}y',z}=N$ for some $z\in W$, then $\gamma_{x'r_{1},r_{1}y',z}=\beta_{x'r_{1},r_{1}y',z}\neq 0$.
By Lemma 4.1, there exists $y''\in W$, such that $y'=x'^{-1}y''$. In this case, we have $z=(w_{J}y'')^{-1}$ and $\beta_{x'r_{1},r_{1}y',(w_{J}y'')^{-1}}=\gamma_{x'r_{1},r_{1}y',(w_{J}y'')^{-1}}=1$,
so precisely one of $deg\ ((v_{r_{1}}-v_{r_{1}}^{-1})T_{x'r_{1}}T_{y'})$ and $deg\ (T_{x'}T_{y'})$ is $N$. Similarly, we know precisely one of $deg\ ((v_{r_{2}}-v_{r_{2}}^{-1})T_{x'r_{2}}T_{y'})$ and $deg\ (T_{x'}T_{y'})$ is $N$.

Now we claim that precisely one of $deg\ ((v_{r_{1}}-v_{r_{1}}^{-1})T_{x'r_{1}}T_{y'})$, $deg\ ((v_{r_{2}}-v_{r_{2}}^{-1})T_{x'r_{2}}T_{y'})$ and $deg\ (T_{x'}T_{y'})$ is $N$. Otherwise, we must have
$deg\ ((v_{r_{1}}-v_{r_{1}}^{-1})T_{x'r_{1}}T_{y'})=deg\ ((v_{r_{2}}-v_{r_{2}}^{-1})T_{x'r_{2}}T_{y'})=N$ and $deg\ (T_{x'}T_{y'})< N$. Thus $\mathcal{L}(x')\neq J$.
By Lemma 4.4, we have $\mathcal{L}(x'r_{1})\neq J$ or $\mathcal{L}(x'r_{2})\neq J$. Since $x'w_{J'}\underset{R}{\leq} x'r_{1}$, $x'w_{J'}\underset{R}{\leq} x'r_{2}$, we get $x'r_{1}\notin c_{0}$ or $x'r_{2}\notin c_{0}$ by $P10'$, contradict to Corollary 2.3(1).

Summarizing the arguments above, we know
$$\sum_{l(q)\leq 1}f_{w_{J'},w_{J'},q}T_{x'q}T_{y'}$$
contributes $\delta t_{w_{J}xy}$ for $t_{x}t_{y}$. If there exists $y_{1}\in W$ such that $y=x^{-1}\cdot y_{1}$, then $\delta=1$, otherwise $\delta=0$.

We have completed the proof of Proposition 4.2.\hfill $\square$


\medskip

\begin{thebibliography}{Lam86}


\bibitem[AF]{bjorner}
A.Bjorner, F.Brenti:
\emph{Combinatorics of Coxeter groups},
Springer (2005).

\bibitem[Bre]{bremke}
K.Bremke:
\emph{On generalized cells in affine Weyl groups},
J.Algebra {\bf 191} (1997) 149--173.

\bibitem[Gao]{gao}
J.Gao:
\emph{The boundness of weighted Coxeter groups of rank 3},
arXiv:1607.02286v10, 2019.

\bibitem[H]{humphreys}
J.E.Humphreys:
\emph{Reflection groups and Coxeter groups},
Cambridge Studies in Advanced Mathematics {\bf 29} (1992).

\bibitem[KL]{kazhdan}
D.Kazhdan, G.Lusztig:
\emph{Representations of Coxeter groups and Hecke algebras},
Invent. Math. {\bf 53} (1979) 165--184.

\bibitem[L1]{lusztig}
G.Lusztig:
\emph{Cells in affine Weyl groups},
Algebraic groups and related topics, Adv. Stud. Pure Math. {\bf 6} (1985) 255--287.

\bibitem[L2]{lusztig}
G.Lusztig:
\emph{Cells in affine Weyl groups, \uppercase\expandafter{\romannumeral2}},
J.Algebra {\bf 109} (1987) 536--548.

\bibitem[L3]{lusztig}
G.Lusztig:
\emph{Hecke algebras with unequal parameters},
CRM monograph series {\bf 18} (2003).

\bibitem[SY]{shi}
J.Shi, G.Yang:
\emph{The boundness of the weighted Coxeter group with complete graph},
Proc. Amer. Math. Soc. {\bf 144} (2016) 4573--4581.

\bibitem[Xi1]{xi}
N.Xi:
\emph{Representations of affine Hecke algebras},
Lect. Notes Math., Springer {\bf 1587} (1994).

\bibitem[Xi2]{xi}
N.Xi:
\emph{Lusztig's a-function for Coxeter groups with complete graphs},
Bull. Inst. Math. Acad. Sinica {\bf 7 (1)} (2012) 71--90.

\bibitem[Xie1]{xie}
X.Xie:
\emph{The based ring of the lowest generalized two-sided cell of an extended affine Weyl group},
J.Algebra {\bf 477} (2017) 1--28.

\bibitem[Xie2]{xie}
X.Xie:
\emph{The lowest two-sided cell of a Coxeter group with complete graph},
J.Algebra {\bf 489} (2017) 38--58.

\bibitem[Zhou]{zhou}
P.Zhou:
\emph{Lusztig's a-function for Coxeter groups of rank 3},
J.Algebra {\bf 384} (2013) 169--193.


\end{thebibliography}
\end{document}